\documentclass[twocolumn]{autart}
\usepackage{cite}
\usepackage[cmex10]{amsmath}
\usepackage{graphicx}
\usepackage{subfigure}
\usepackage{amsfonts}
\usepackage{enumerate}

\begin{document}

\begin{frontmatter}

\title{Stochastic Source Seeking with Forward and Angular Velocity Regulation \thanksref{footnoteinfo}}
\thanks[footnoteinfo]{The material in this paper was not presented at any conference.}

\author[Beijing]{Jinbiao Lin}\ead{linjb11@mails.tsinghua.edu.cn},
\author[Beijing]{Shiji~Song\corauthref{cor}}\ead{shijis@tsinghua.edu.cn},
\author[Beijing]{Keyou~You}\ead{youky@tsinghua.edu.cn},
\author[SanDiego]{Miroslav Krstic}\ead{krstic@ucsd.edu}
\address[Beijing]{Department of Automation and Tsinghua National Laboratory for Information Science and Technology, Tsinghua University, Beijing, 100084, China}
\address[SanDiego]{Department of Mechanical and Aerospace Engineering, University of California, San Diego, La Jolla, CA 92093-0411, USA}
\corauth[cor]{Corresponding author. Tel. +86010-62782721. Fax +86010-62782721.}


\begin{abstract}                          
This paper studies a stochastic extremum seeking method to steer a nonholonomic vehicle to the unknown source of a static spatially distributed filed in a plane. The key challenge lies in the lack of vehicle's position information and the distribution of the scalar field. Different from the existing stochastic strategy that keeps the forward velocity constant and controls only the angular velocity, we design a stochastic extremum seeking controller to regulate both forward and angular velocities simultaneously in this work. Thus, the vehicle decelerates near the source and stays within a small area as if it comes to a full stop, which solves the overshoot problem in the constant forward velocity case. We use the stochastic averaging theory to prove the local exponential convergence, both almost surely and in probability, to a small neighborhood near the source for elliptical level sets. Finally, simulations are included to illustrate the theoretical results.
\end{abstract}

\end{frontmatter}

\section{Introduction}
{\it Source seeking} is a problem of steering single or multiple autonomous agents to seek the source of an unknown scalar field, which may be thermal, electromagnetic, acoustic, or the concentration of a chemical agent. Source seeking is of interest in many areas, such as environmental studies, explosive detection, localizing the sources of hazardous chemicals leakage or pollutants, etc. There are a diversity of approaches to source seeking problem. Motivated by biological chemotactic and anemotactic behaviors, behavior-based adaptive mission planners are proposed to trace a chemical plume to locate its source \cite{russell2004robotic,li2006moth}. In an alternative approach, mathematical programming methods, such as gradient descent method, are adopted to address this problem \cite{Porat1996,tee2001solving,Ogren2004,khong2014multi}. Besides, source-likelihood mapping methods are studied in \cite{jakuba2007stochastic,pang2006chemical}.

In this work, we consider steering a {\it single nonholonomic} vehicle to locate a {\it static} source which creates a continuous signal map {\it in a plane}. Recently there is a growing interest in the study of locating such a source {\it without position information} \cite{matveev2011navigation,azuma2012stochastic,Cochran2009}. The lack of position information is taken account for vehicles operated in environments where their position information is unavailable or costly. Obviously, this constraint, along with the nonholonomic constraint of the vehicle kinematics, renders the guidance of the vehicle interesting and challenging.

{\it Extremum seeking} (ES) is a model free optimization method for dynamical system with limited knowledge\cite{krstic2000stability}. It has been proved to an effective method for nonholonomic source seeking problems without position information. In \cite{Zhang2007} and \cite{Cochran2009} ES was applied to tune the forward or angular velocity of the vehicle to locate the source. In \cite{Ghods2010} Ghods and Krstic regulated both velocities to control the vehicle to stop near the source. While the above works focus on the 2D vehicles, Lin et al. \cite{Lin2015} considered the more complicated 3D case. Different from the above perturbation-based ES methods, a novel regulator without injecting any perturbation is proposed based on Lie bracket approximation in \cite{scheinker2014extremum,durr2016}.

Motivated by the chemotactic behavior of the bacterium Escherichia coli ({\it E. coli}) \cite{berg2008coli}, Liu and Krstic applied the stochastic averaging technique to ES \cite{liu2010stochastic}. In \cite{Liu2010}, the stochastic ES algorithm is applied to the source seeking problem. Thus, the seeker can successfully locate the source but with an unpredictable, ``nearly random" trajectory. This feature would be useful when the seeker itself is pursued by another hostile pursuer. Whereas in \cite{Liu2010} the forward velocity is chosen to be constant, which results in complicated asymptotic behaviors. Particularly, the vehicle cannot settle when it approaches close to the source. Instead it exhibits certain overshoots and finally revolves around the source. A small constant forward velocity may improve the asymptotic performance, but it will decrease the convergence rate.

In order to improve the asymptotic performance of the vehicle, we applied the stochastic ES to tune both forward and angular velocities simultaneously, which is different from \cite{Liu2010}. Note that the deterministic case has been discussed in \cite{Ghods2010}, and in this work we focus on designing a stochastic excitation to modulate the velocities. Under a tunable forward velocity we can slow down the vehicle around the source to approach the source. In addition, the undesired overshoots are eliminated due to a tunable forward velocity. We adopt the stochastic averaging theory to establish the local exponential convergence, both almost surely and in probability, to a small neighborhood near the source, for signal fields with elliptical level sets. Note that in \cite{Cochran2009,Liu2010,Ghods2010} only the stability for circular level sets was proved, and in this work the stability for elliptical level sets is firstly proved\footnote{In \cite{scheinker2014extremum,durr2016}, the stability for a general signal map was proved under a different strategy, but their approach is inefficient in the convergence towards the source.}. 

It should be mentioned that besides stochastic ES there are several methods to address stochastic source seeking without position information. In \cite{azuma2012stochastic} Azuma et al. adopted the stochastic approximation technique to solve this problem by sequentially generating waypoints which converges to the source. Their method can work for a switching signal field under the assumption that the robot can move to any point in the body fixed coordinate frame. Apparently the controller is discontinuous and it requires accurate clocks to decide whether to enter next step. Another representative method is proposed in \cite{mesquita2008optimotaxis}, where induces a swarm of autonomous vehicles to to perform a biased random walk and finally achieves higher vehicle densities near the maximum. The convergence of the agents probability density to a specified function of the spatial profile of the measured signal was demonstrated. Compared with these methods, the advantage of the stochastic ES is that we can use a single vehicle with a simple continuous controller to locate the source. What's more, the exponential convergence (in probability and almost surely) to to a small attractor near the source can be established. However, the convergence results of stochastic ES are only for static quadratic signal maps while in \cite{azuma2012stochastic,mesquita2008optimotaxis} the signal map can be more general.

The rest of the paper is organized as follows. In Section \ref{sec:problem} we describe the nonholonomic source seeking problem and propose the stochastic ES scheme. In Section \ref{sec:stability} we prove the local exponential convergence for signal fields with elliptical level sets. We first derive an average system to approximate the original system, then we consider the local stability under different bias forward velocities. After that we discuss the result for circular level sets as a special case. In Section \ref{sec:simulation} we include simulation results to illustrate the effectiveness of the control scheme.

\section{Problem Description and Control Scheme} \label{sec:problem}
In this section we firstly describe the vehicle model and formulate the source seeking problem. Then we propose a stochastic ES scheme to adjust the forward and angular velocities of the vehicle.
\subsection{Problem Description}
\begin{figure}[!t]
  \centering
  \includegraphics[height=4cm]{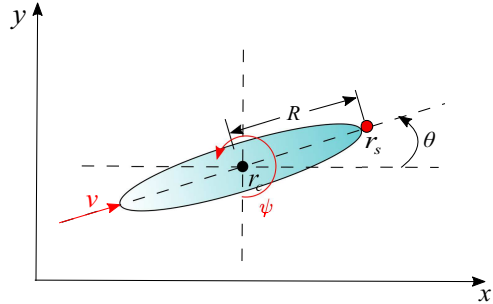}
  \caption{Geometric interpretation of vehicle model.}
  \label{figModel}
\end{figure}
Similar to \cite{Cochran2009}, we consider an autonomous vehicle modeled as a nonholonomic unicycle,  see Fig.~\ref{figModel} for illustration. The heading angle is defined by $\theta$, and the position of the vehicle center is defined by $r_c$. A sensor is mounted at the front end $r_s$, a distance $R$ away from the vehicle center $r_c$. The vehicle has actuators which are used to impart the forward velocity $v$ and the angular velocity $\psi$.
The kinematic equations of motion for the vehicle center and the sensor are
\begin{align}
{{\dot r}_c} &= v{{\rm e}^{j\theta }}, \label{eq:motion1}\\
 \dot \theta  &= \psi, \\
 {r_s} &= {r_c} + R{{\rm e}^{j\theta }},\label{eq:motion3}
\end{align}
where $r_c$ and $r_s$ are written as complex variables.

The task of vehicle is to seek a static source that emits a spatially distributed signal in a plane. We denote the signal strength at the location $r$ by $J = f\left( r \right)$ and make the following assumption.
\begin{assum}\label{assumption1}
The signal strength $J$ decays away from the source and achieves its isolated local maximum $f^* = f(r^*)$ at the source location $r^*$. What's more, the distribution $f\left( r \right)$ is twice continuously differentiable and
\begin{align}
&\nabla f(r^*)=0,\\
&\nabla^2 f(r^*)\mbox{\;is negative definite},
\end{align}
where $\nabla f(r^*)$ and $\nabla^2 f(r^*)$ denote the gradient and Hessian of $f$ at $r^*$ respectively.
\end{assum}
Assumption \ref{assumption1} is a natural extension of Assumption 2.3 in \cite{krstic2000stability}. Under this assumption, we can approximate the signal distribution by a quadratic map when studying the local convergence. Without loss of generality, we assume the quadratic map takes the form
\begin{equation}
J = {f^*} - q_x( x_s  - {x^*} )^2 - q_y( y_s  - {y^*} )^2, \label{eq:eplliticalMap1}
\end{equation}
where $r_s=[x_s,y_s]^T$, $r^*=[x^*,y^*]^T$, and $q_x$ and $q_y$ are unknown positive constants.

The objective of this work is to design a control scheme to navigate the vehicle to the unknown source. The signal strength $J$ can be measured by the sensor, but the position information is unavailable. The explicit form of the signal field, such as the shape of $f$ and the position of $r^*$, is also unknown. Note that the traditional gradient searching strategy is not suitable for the problem due to the lack of position information.

\subsection{Control Scheme}
\begin{figure}[!t]
  \centering
  \includegraphics[width=0.8\hsize]{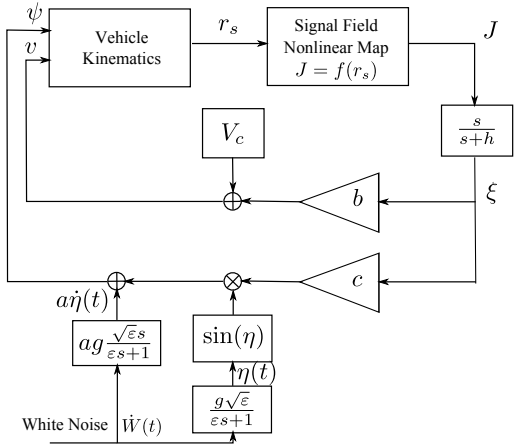}
  \caption{Block diagram of stochastic source seeking via tuning of the forward and angular velocities of the vehicle.}
  \label{figController}
\end{figure}
We employ the stochastic ES method to tune the angular velocity $\psi$ directly and the forward velocity $v$ indirectly. The control scheme is depicted in Fig.~\ref{figController}. The control laws are given by
\begin{align}
v &= {V_c} + b\xi, \label{eq:law1}\\
\psi  &= a\dot \eta + c\xi \sin (\eta ), \label{eq:law2}\\
\xi  &= \frac{s}{{s + h}}[J], \label{eq:law3}\\
\eta  &= \frac{{g\sqrt \varepsilon  }}{{\varepsilon s + 1}}[{\dot W}], \label{eq:law4}
\end{align}
where the parameters $a$, $g$, $\varepsilon$, $b$, $c$, $h$ and $V_c$ are positive and will affect the performance of the approach, $J$ is the sensor reading, and $W(t)$ is a standard Brownian motion defined in a complete probability space $(\Omega, \mathcal{F}, P)$ with the sample space $\Omega$, the $\sigma$-field $\mathcal{F}$, and the probability measure $P$. Here the colored noise $\eta$ is used as a stochastic perturbation in ES.

In our control scheme, the angular velocity $\psi$ is tuned according to the idea of the stochastic ES tuning law \cite{liu2010stochastic}. The perturbation term $a\dot \eta$ is added to persistently excite the system while the corresponding demodulation term $\sin (\eta)$ is used to estimate the gradient of the nonlinear map $f$. Different from the deterministic case which uses a sinusoidal perturbation \cite{Cochran2009,Ghods2010}, the stochastic perturbation results in a partly random trajectory. The forward velocity $v$ is designed to be positively correlated to $\xi$, since $\xi$ describes the variation of the sensor reading $J$ in some sense. As a result, the vehicle would speed up when approaching the source, and slow down when deviating from the source.

It is worthy mentioning that our control scheme is different from the one in \cite{Liu2010}, where a vehicle with a constant forward velocity is considered. Employing the basic stochastic ES method in \cite{liu2010stochastic}, the vehicle with a constant forward velocity cannot settle even if it has reached the source. In addition, it easily overshoots the source and has to turn around. This process may repeat for a while before the vehicle finally revolves around the source. In \cite{Liu2010} a nonlinear damping item is added to tune the angular velocity to improve the performance and achieve exponential stability. In this work, we tune the forward velocity along with the angular velocity. Intuitively this is a better way to control the vehicle, as we are able to smartly adjust the vehicle to speed up or slow down depending on different circumstances. Thus we can apply the basic stochastic ES control law to tuning the angular velocity directly without employing the nonlinear damping.

\section{Stability Analysis} \label{sec:stability}
The dynamics of the closed-loop system is intricate on the account of nonlinearities of the vehicle model and the signal map and the existence of the stochastic perturbation. We adopt the stochastic averaging theory in \cite{liu2010stochastic} to prove the local exponential convergence for elliptical level sets. In Section \ref{subsec:average} we derive an average system to approximate the closed-loop system. In Section \ref{subsec:elliptical} we prove that the vehicle converges, almost surely and in probability, to an attractor near the source under a small bias forward velocity. Section \ref{subsec:circular} we consider a special case where the signal distribution is circular. In this case, the local exponential convergence can be established no matter the bias forward velocity is small or large.

\subsection{Average System for Elliptical Level Sets} \label{subsec:average}
We firstly rewrite the elliptical signal map (\ref{eq:eplliticalMap1}) as
\begin{align*}
J & = {f^*} - ({q_r}+2{q_p})( x_s  - {x^*} )^2 - ({q_r}-2{q_p})( y_s  - {y^*} )^2\\
 &= {f^*} - {q_r}{\left| {{r_s} - {r^*}} \right|^2} - {q_p}\left( {{{({r_s} - {r^*})}^2} + {{(\overline {{r_s} - {r^*}} )}^2}} \right),
\end{align*}
where $q_r$ and $q_p$ are unknown and $q_r>2|q_p|\ge0$.

Before analyzing the stability of the closed-loop system, we define an output error variable $e_\xi = \frac{h}{{s + h}}[J] - {f^*}$ to express the output of the washout filter as $$\xi  = \frac{s}{{s + h}}[J] = J - \frac{h}{{s + h}}[J] = J - {f^*} - e_\xi.$$ Thus, we obtain $\dot e_\xi = h\xi$.

By inserting the control laws (\ref{eq:law1})-(\ref{eq:law4}) into the system (\ref{eq:motion1})-(\ref{eq:motion3}) and expressing $\dot\eta$ as
\begin{equation*}
\dot\eta = \frac{{g\sqrt \varepsilon s }}{{\varepsilon s + 1}}[{\dot W}]
=\frac{1}{\sqrt \varepsilon} \frac{{g\varepsilon s + g - g}}{{\varepsilon s + 1}}[{\dot W}]
=\frac{g}{\sqrt \varepsilon}{\dot W} - \frac{1}{\varepsilon} \eta,
\end{equation*}
the closed-loop system is written as
\begin{align}
{\rm d}r_c &= ({V_c} + b\xi ){\rm e}^{j\theta}{\rm d}t, \label{eq:sys1}\\
{\rm d} \theta  &=  -\frac{a}{\varepsilon}\eta{\rm d}t + c\xi \sin (\eta ){\rm d}t + \frac{ag}{\sqrt \varepsilon}{\rm d}W,\\
{\rm d} e_\xi &= h\xi,\\
{\rm d} \eta &= - \frac{1}{\varepsilon} \eta {\rm d} t + \frac{g}{\sqrt \varepsilon}{{\rm d} W}, \label{eq:sys4}\\
\xi & = J - {f^*} - e_\xi,\\
{r_s} &= {r_c} + R{\rm e}^{j\theta}. \label{eq:sys6}
\end{align}

To analyze the closed-loop system, we firstly re-express it by variable transformation. Then we redefine $r_c$ in its polar coordinates for the convenience of the calculation of the equilibria. To this end, we start by defining shifted variables
$${\hat r}_c = r_c - r^*, \quad \hat \theta = \theta  - a\eta.$$
The dynamics of the shifted system is given by
\begin{align*}
{\rm d}\hat r_c &= ({V_c} + b\xi ){{\rm e}^{j(\hat \theta  + a\eta)}},\\
{\rm d}\hat \theta  &=  c\xi \sin ( \eta ){\rm d}t,\\
{\rm d} e_\xi &= h\xi,\\
\xi  &= - e_\xi - {q_r}{\left| {\hat r}_c + R{e^{j\theta}} \right|^2} \\
&\quad - {q_p}\left( ({{\hat r}_c + R{{\rm e}^{j\theta}})^2} + {{(\overline{ {\hat r}_c + R{{\rm e}^{j\theta}} })}^2} \right),\\
{\rm d} \eta &= - \frac{1}{\varepsilon} \eta {\rm d} t + \frac{g}{\sqrt \varepsilon}{{\rm d} W}.
\end{align*}

By (\ref{eq:sys4}) and the definition of Ito stochastic differential equation, we obtain $\eta (t) = \eta (0) - \int_0^t {\frac{1}{\varepsilon }\eta (\tau){\rm d}\tau}  + \int_0^t {\frac{g}{{\sqrt \varepsilon  }}{\rm d}W(\tau)}$. Then we have $\eta (\varepsilon t) = \eta (0) - \int_0^t {\eta (\varepsilon u){\rm d}u}  + \int_0^t {\frac{g}{{\sqrt \varepsilon  }}{\rm d}W(\varepsilon u)}$. By defining $\chi (t) = \eta (\varepsilon t)$ and $B(t) = \frac{1}{{\sqrt \varepsilon  }}W(\varepsilon t)$, we have $ {\rm d}\chi (t) =  - \chi (t){\rm d}t + g{\rm d}B(t)$, where $B(t)$ is a standard Brownian motion and the process $\chi(t)$ is an Ornstein--Uhlenbeck (OU) process which is ergodic with invariant distribution $\mu ({\rm d}y) = \frac{1}{{\sqrt \pi  g}}{{\rm e}^{ - \frac{{{y^2}}}{{{g^2}}}{\rm d}y}}$.

We redefine ${\hat r}_c$ by its polar coordinates
\begin{align*}
- {{\hat r}_c} &= \left| {\hat r_c} \right|{{\rm e}^{j{\theta ^*}}} = {\tilde r_c}{{\rm e}^{j{\theta ^*}}}, \\
\theta^* &= \arg ( - {\hat r_c}) = \arg (r^*-r_c),
\end{align*}
where ${\tilde r}_c$ is the distance between the vehicle center and the source, ${\theta ^*}$ represents the heading angle from the vehicle center towards the source. We also define $\tilde e_\xi  = e_\xi + q_rR^2$ for convenience. Using these new definitions, $\xi$ is expressed as
\begin{align}
\xi  &=  - {\tilde e_\xi} - { 2{q_p}R^2\cos (2\hat \theta  + 2a\chi (t/\varepsilon ) )}  \notag\\
&\quad - {{\tilde r}^2}\left( {{q_r} + 2{q_p}\cos (2{\theta ^*})} \right) \notag\\
&\quad + 2\tilde rR\big( {q_r}\cos (\hat \theta  - {\theta ^*} + a\chi (t/\varepsilon ) ) \notag\\
&\quad\quad\quad\quad\;+ 2{q_p}\cos (\hat \theta  + {\theta ^*} + a\chi (t/\varepsilon ) ) \big).\label{eq:nxi}
\end{align}
Now we obtain the following shifted error system
\begin{align}
\frac{{\rm d}{{\tilde r}_c}}{{\rm d}t} &=  - ({V_c} + b\xi )\cos (\hat\theta - \theta^*  + a\chi (t/\varepsilon )), \label{eq:eplliticalSys1}\\
\frac{{\rm d} \theta^* }{{\rm d}t} &= - \frac{{V_c} + b\xi }{{{{\tilde r}_c}}}\sin(\hat\theta - \theta^*  + a\chi (t/\varepsilon )),\\
\frac{{\rm d}\hat \theta }{{\rm d}t} &= c\xi \sin (\chi (t/\varepsilon )),\\
\frac{{\rm d}{{\tilde e}_\xi }}{{\rm d}t} &= h\xi,\\
{\rm d}\chi (t) &=  - \chi (t){\rm d}t + g{\rm d}B(t). \label{eq:eplliticalSys5}
\end{align}
According to the stochastic averaging theory \cite{liu2010stochastic}, the error system can be approximated by an average system, which is given by
\begin{align}
\frac{{{\rm d}{{\tilde r}_c}^{\rm av}}}{{{\rm d}t}} &= \left( {b\phi_0 - {V_c}} \right) \cos ({\hat \theta}^{\rm av}  - {\theta ^*}^{\rm av}){I_1}(a,g) \notag\\
&\quad + b{q_p}{R^2}\phi_1 - b{{\tilde r}_c}^{\rm av} R\phi_2, \label{eq:eplliticalEASys1}\\
\frac{{{\rm d}{\theta ^*}^{\rm av}}}{{{\rm d}t}} &= \frac{b\phi_0 - {V_c}}{{{{{\tilde r}_c}^{\rm av}}}}\sin({\hat \theta}^{\rm av} - {\theta ^*}^{\rm av}){I_1}(a,g) \notag\\
&\quad + \frac{{b{q_p}{R^2}}}{{{\tilde r}_c}^{\rm av}}\phi_3 - bR\phi_4,\\
\frac{{{\rm d}{\hat \theta}^{\rm av} }}{{{\rm d}t}} &= 2c{q_p}{R^2}\sin(2{\hat \theta}^{\rm av} ){I_2}(2a,g) \notag\\
&\quad - 2c\tilde rR \phi_5 {I_2}(a,g),
\end{align}
\begin{align}
\frac{{{\rm d}{{\tilde e}_\xi }}}{{{\rm d}t}} &=  - 2h{q_p}{R^2}\cos (2{\hat \theta}^{\rm av} ){I_1}(2a,g) - h\phi_0 \notag\\
&\quad + 2h{{\tilde r}_c}^{\rm av} R\phi_6 {I_1}(a,g), \label{eq:eplliticalEASys4}
\end{align}
where ${I_1}(a,g) = \int\limits_\Re  {\cos (ay)\mu (dy)}  = {e^{ - \frac{{{a^2}{g^2}}}{4}}} $, ${I_2}(a,g) $$= \int\limits_\Re  {\sin (ay)\sin(y)\mu ({\rm d}y)} = \frac{1}{2}\left[ {{e^{ - \frac{{{{(a - 1)}^2}{g^2}}}{4}}} - {e^{ - \frac{{{{(a + 1)}^2}{g^2}}}{4}}}} \right]$ and
\begin{align*}
\phi_0 &= {{\tilde e}_\xi }^{\rm av} + ({{\tilde r}_c}^{\rm av})^2({q_r} + 2{q_p}\cos (2{\theta ^*}^{\rm av})),\\
\phi_1 &= \cos (3{\hat \theta}^{\rm av}  - {\theta ^*}^{\rm av}){I_1}(3a,g)\\
&\quad + \cos ({\hat \theta}^{\rm av}  + {\theta ^*}^{\rm av}){I_1}(a,g),\\
\phi_2 &= {q_r}\cos (2{\hat \theta}^{\rm av}  - 2{\theta ^*}^{\rm av}){I_1}(2a,g) + 2{q_p}\cos (2{\theta ^*}^{\rm av}) \\
&\quad+ 2{q_p}\cos (2{\hat \theta}^{\rm av} ){I_1}(2a,g) + {q_r},\\
\phi_3 &= \sin(3{\hat \theta}^{\rm av}  - {\theta ^*}^{\rm av}){I_1}(3a,g) - \sin({\hat \theta}^{\rm av}  + {\theta ^*}^{\rm av}){I_1}(a,g),\\
\phi_4 &= {q_r}\sin (2{\hat \theta}^{\rm av}  - 2{\theta ^*}^{\rm av}){I_1}(2a,g) - 2{q_p}\sin(2{\theta ^*}^{\rm av}) \notag\\
&\quad+ 2{q_p}\sin (2{\hat \theta}^{\rm av} ){I_1}(2a,g),\\
\phi_5 &= {q_r}\sin({\hat \theta}^{\rm av}  - {\theta ^*}^{\rm av}) + 2{q_p}\sin({\hat \theta}^{\rm av}  + {\theta ^*}^{\rm av}),\\
\phi_6 &= {q_r}\cos({\hat \theta}^{\rm av}  - {\theta ^*}^{\rm av}) + 2{q_p}\cos({\hat \theta}^{\rm av}  + {\theta ^*}^{\rm av}).
\end{align*}

The average error system has eight equilibria as follows\footnote{We have implicitly assumed ${{\theta^* }^{{\rm{ave}}}} \in (-\pi,\pi]$ and ${{\hat \theta }^{{\rm{ave}}}} \in (-\pi,\pi]$ to exclude repetitive equilibria.}
\begin{align}
{\bf eq_1} &= \left[ {\rho(q_p)},0,0,{e(q_p)} \right], \label{eq:neq1}\\
{\bf eq_2} &= \left[ {\rho(q_p)},\pi,\pi,{e(q_p)} \right], \label{eq:neq2}\\
{\bf eq_3} &= \left[ {\rho(-q_p)},\pi/2,\pi/2,{e(-q_p)} \right], \label{eq:neq3}\\
{\bf eq_4} &= \left[ {\rho(-q_p)},-\pi/2,-\pi/2,{e(-q_p)} \right], \label{eq:neq4}\\
{\bf eq_5} &= \left[ {-\rho(q_p)},\pi,0,{e(q_p)} \right], \label{eq:neq5}\\
{\bf eq_6} &= \left[ {-\rho(q_p)},0,\pi,{e(q_p)} \right], \label{eq:neq6}\\
{\bf eq_7} &= \left[ {-\rho(-q_p)},-\pi/2,\pi/2,{e(-q_p)} \right], \label{eq:neq7}\\
{\bf eq_8} &= \left[ {-\rho(-q_p)},\pi/2,-\pi/2,{e(-q_p)} \right], \label{eq:neq8}
\end{align}
where ${\bf eq_i}$ is of the form $[ {\tilde r_c}^{\rm av},{\theta ^*}^{\rm av},{\hat \theta}^{\rm av},{\tilde e_\xi}^{{\rm av}} ]$ and
\begin{align*}
\rho ({q_p}) & \buildrel \Delta \over = \frac{{ - {V_c}{I_1}(a,g) + b{q_p}{R^2}{\gamma _2}}}{{bR({q_r} + 2{q_p}){\gamma _1}}},\\
e({q_p}) & \buildrel \Delta \over = 2R({q_r} + 2{q_p}){I_1}(a,g)\rho ({q_p}) - ({q_r} + 2{q_p}){\rho ^2}({q_p})\\
&\quad - 2{q_p}{R^2}{I_1}(2a,g),\\
{\gamma _1} &= 1 + {I_1}(2a,g) - 2I_1^2(a,g),\\
{\gamma _2} &= {I_1}(3a,g) + {I_1}(a,g) - 2{I_1}(2a,g){I_1}(a,g).
\end{align*}
Note that with positive $a$ and $g$, we have $I_1(a,g)>0$, $I_2(a,g)>0$, $\gamma_1>0$ and $\gamma_2>0$.

Each of the equilibria (\ref{eq:neq1})-(\ref{eq:neq8}) represents an attractor around the source. The value of ${\tilde r_c^{\rm{av}}}$ should be real and positive as it represents the average distance between the vehicle center and the source. Note the difference between $\hat \theta^{{\rm{av}}}$ and ${\theta^*}^{{\rm{av}}}$ is either $0$ or $\pi$, which indicates the average heading of the vehicle points either directly towards or away from the source.

\subsection{Stability for Elliptical Level Sets} \label{subsec:elliptical}
Before we declare the stability, we define an index variable $\iota$ as
$$
\iota  = \left\{ {\begin{array}{rcl}
1, & &\text{if\;} {q_p} \gamma_3 < 0 \text{\;and\;} \rho(q_p)>0\\
3, & &\text{if\;} {q_p} \gamma_3 > 0 \text{\;and\;} \rho(-q_p)>0\\
5, & &\text{if\;} {q_p} \gamma_3 < 0 \text{\;and\;} \rho(q_p)<0\\
7, & &\text{if\;} {q_p} \gamma_3 > 0 \text{\;and\;} \rho(-q_p)<0
\end{array}} \right. ,$$
where
\begin{align*}
\gamma_3 &= \left( {{I_1}(3a,g) - {I_1}(a,g)} \right){I_2}(a,g) \\
&\quad + \left(1 - {I_1}(2a,g) \right){I_2}(2a,g).
\end{align*}
\begin{thm} \label{pro3}
Consider the system (\ref{eq:eplliticalMap1}), (\ref{eq:sys1})-(\ref{eq:sys6}) with positive parameters $a$, $g$, $b$, $c$, $h$, and $\varepsilon \in (0,\varepsilon_0)$. The parameters $a$, $g$, $b$, $c$, $h$, $V_c$ are chosen such that either
\begin{gather}
{q_p} \gamma_3 < 0, \label{eq:ConvergeConditionE1} \\
V_c^l(q_p)<V_c<V_c^u(q_p), \text{\quad and \quad} V_c \neq V_{io}, \label{eq:ConvergeConditionE2}
\end{gather}
or
\begin{gather}
{q_p} \gamma_3 > 0, \label{eq:ConvergeConditionE3} \\
V_c^l(-q_p)<V_c<V_c^u(-q_p), \text{\quad and \quad} V_c \neq V_{io},  \label{eq:ConvergeConditionE4}
\end{gather}
where
\begin{align*}
V_c^l(q_p) & \buildrel \Delta \over = \frac{-1}{{2I_1^2(a,g)}}\big( b{R^2}\left( {{q_r} + 2{q_p}} \right){\gamma _1}\left( {1 + {I_1}(2a,g)} \right) \\
&\quad\quad\quad\quad\quad\quad + hR{\gamma _1} - 2b{R^2}{q_p}{\gamma _2}{I_1}(a,g) \big), \\
V_c^u(q_p) & \buildrel \Delta \over = \frac{ {b^2}R({q_r} - 2{q_p}){\gamma _4} + 2bc{q_p}{R^2}\gamma_5 }{{2c{I_1}(a,g){I_2}(a,g)}}, \\
V_{io} & \buildrel \Delta \over = -{\rm sgn}(\gamma_3q_p){b{q_p}{R^2}{\gamma _2}}/I_1(a,g),\\
\gamma_4 &= \gamma_1\left( {1 - {I_1}(2a,g)} \right), \\
\gamma_5 &= {\gamma _2}{I_2}(a,g) - 2{\gamma _1}{I_2}(2a,g).
\end{align*}
If the initial conditions ${r_c}(0),\theta (0),e_\xi(0)$ are such that either $\left| {\bf\Xi}(0) -{\bf eq_\iota} \right|$ or $\left| {\bf\Xi}(0) -{\bf eq_{\iota+1}} \right|$ is sufficiently small, where
\[{\bf\Xi}(t) = \left[|{{r_c}(t) - r^*}|, \arg(r^*-r_c(t)), \theta (t), e_\xi(t) - {q_r}{R^2}\right],\]
then there exist constants $C_0$, $\gamma_0>0$, $T(\varepsilon):(0,\varepsilon_0)\to \mathbb{N}$ such that for any $\delta>0$, the trajectory of the vehicle center $r_c(t)$ satisfies the following properties,
\begin{align*}
\mathop {\lim }\limits_{\varepsilon  \to \infty } &\inf \big\{ t \ge 0: \left| {| {{r_c}(t) - r^*} | - \tilde r_\iota } \right|\\
&\quad\quad\quad\quad\quad\quad > {C_0}{{\rm e}^{ - {\gamma _0}t}} + \delta \big\} = \infty ,\quad \text{a.s.},\\
\mathop {\lim }\limits_{\varepsilon  \to \infty } &P \big\{ \left| {| {{r_c}(t) - r^*} | - \tilde r_\iota } \right| \le {C_0}{{\rm e}^{ - {\gamma _0}t}} + \delta,\\
&\quad \forall t\in[0,T(\varepsilon)] \big\} =1 \text{\quad with} \mathop {\lim }_{\varepsilon  \to \infty }T(\varepsilon)=\infty,
\end{align*}
where $\tilde r_\iota$ denotes the first element of ${\bf eq_\iota}$, the constant $C_0$ is dependent on the initial condition $\left( r_c(0),\theta(0),e(0) \right)$ and on the parameters $a$, $g$, $b$, $c$, $h$, $V_c$, $R$, $q_r$, $q_p$, and the constant $\gamma_0$ is dependent on the parameters $a$, $g$, $b$, $c$, $h$, $V_c$, $R$, $q_r$, $q_p$.
\end{thm}
\begin{pf}
The Jacobians of the equilibria (\ref{eq:neq1})-(\ref{eq:neq8}) are as follows,
\begin{align*}
{A}^{\bf eq1} &= {A}^{\bf eq2} = J_1(q_p), &\quad {A}^{\bf eq3} &= {A}^{\bf eq4} = J_1(-q_p), \\
{A}^{\bf eq5} &= {A}^{\bf eq6} = J_2(q_p), &\quad {A}^{\bf eq7} &= {A}^{\bf eq8} = J_2(-q_p),
\end{align*}
where $J_1(q_p)$ and $J_2(q_p)$ are defined as
\begin{align*}
J_1(q_p) & \buildrel \Delta \over = \begin{bmatrix}
{{a_{11}}({q_p})}&0&0&{{a_{14}}}\\
0&{{a_{22}}({q_p})}&{{a_{23}}({q_p})}&0\\
0&{{a_{32}}({q_p})}&{{a_{33}}({q_p})}&0\\
{{a_{41}}({q_p})}&0&0&{ - h}
\end{bmatrix},\\
J_2(q_p) & \buildrel \Delta \over = \begin{bmatrix}
{{a_{11}}({q_p})}&0&0&{{-a_{14}}}\\
0&{{a_{22}}({q_p})}&{{a_{23}}({q_p})}&0\\
0&{{a_{32}}({q_p})}&{{a_{33}}({q_p})}&0\\
{{-a_{41}}({q_p})}&0&0&{ - h}
\end{bmatrix},
\end{align*}
and the explicit forms of $a_{ij}$ are given in Appendix \ref{appendix:coefficient}.

The Jacobians ${A}^{\bf eq1}$, ${A}^{\bf eq2}$, ${A}^{\bf eq5}$ and ${A}^{\bf eq6}$ have the same characteristic equation, which is given by
\begin{align}
&\left[ {{\lambda ^2} + (h - {a_{11}(q_p)})\lambda  - {a_{11}(q_p)}h - {a_{14}(q_p)}{a_{41}(q_p)}} \right] \notag\\
&\times \big[ {\lambda ^2} - ({a_{22}(q_p)} + {a_{33}(q_p)})\lambda  + {a_{22}(q_p)}{a_{33}(q_p)} \notag\\
&\quad\; - {a_{23}(q_p)}{a_{32}(q_p)} \big] = 0. \label{eq:characteristic3}
\end{align}
To guarantee that all roots of characteristic equation (\ref{eq:characteristic3}) have negative real parts, we need
\begin{align*}
{a_{11}(q_p)} - h &< 0, \\
{a_{11}(q_p)}h + {a_{14}}{a_{41}(q_p)} &< 0, \\
{a_{22}(q_p)} + {a_{33}(q_p)} &< 0, \\
{a_{23}(q_p)}{a_{32}(q_p)} - {a_{22}(q_p)}{a_{33}(q_p)} &< 0.
\end{align*}
All the above requirements are satisfied under conditions (\ref{eq:ConvergeConditionE1}) and (\ref{eq:ConvergeConditionE2}), which implies the Jacobians ${A}^{\bf eq1}$, ${A}^{\bf eq2}$, ${A}^{\bf eq5}$ and ${A}^{\bf eq6}$ are Hurwitz. Hence, equilibria ${\bf eq_1}$, ${\bf eq_2}$, ${\bf eq_5}$ and ${\bf eq_6}$ are exponentially stable. Similarly we can prove that equilibria ${\bf eq_3}$, ${\bf eq_4}$, ${\bf eq_7}$ and ${\bf eq_8}$ are exponentially stable
under conditions (\ref{eq:ConvergeConditionE3}) and (\ref{eq:ConvergeConditionE4}).
By Theorem 2 in \cite{liu2010stochastic}, there exist constants $c_0^{(i)}>0$, $r_0^{(i)}>0$, $\gamma_0^{(i)}>0$ and a function $T^{(i)}(\varepsilon):(0,\varepsilon_0)\to \mathbb{N}$, $i=1,\cdots,8$, such that for any $\delta>0$ and any initial condition $| {\Lambda _\varepsilon ^{(i)}(0)} | < r_0^{(i)}$,
\begin{align*}
\mathop {\lim }\limits_{\varepsilon  \to \infty } &\inf \big\{ t \ge 0: | {{\bf\Xi} _\varepsilon ^{(i)}(t)} | \\
&\quad\quad > {c_0^{(i)}}| {{\bf\Xi} _\varepsilon ^{(i)}(0)}| {\rm e}^{-\gamma_0^{(i)}t} + \delta \big\} = \infty ,\text{\quad a.s.},\\
\mathop {\lim }\limits_{\varepsilon  \to \infty } &P \big\{ | {{\bf\Xi} _\varepsilon ^{(i)}(t)} | \le {c_0^{(i)}}| {{\bf\Xi} _\varepsilon ^{(i)}(0)}| {\rm e}^{-\gamma_0^{(i)}t} + \delta,\notag\\
&\quad \forall t\in[0,T^{(i)}(\varepsilon)] \big\} =1 \text{\quad with} \mathop {\lim }_{\varepsilon  \to \infty }T^{(i)}(\varepsilon)=\infty,
\end{align*}
where ${{\bf\Xi} _\varepsilon ^{(i)}(t)} = \left| {\bf\Xi}(t) -{\bf eq_i} \right|$. With the fact $|\tilde r_c(t)- \tilde r_\iota| < |{{\bf\Xi} _\varepsilon ^{(\iota)}(t)}|$, we obtain
\begin{align*}
\mathop {\lim }\limits_{\varepsilon  \to \infty } &\inf \big\{ t \ge 0: | \tilde r_c(t)-\tilde r_\iota | \\
&\quad\quad\quad\quad\quad\quad > {C_0^{(\iota)}} {\rm e}^{-\gamma_0^{(\iota)}t} + \delta \big\} = \infty ,\quad \text{a.s.},\\
\mathop {\lim }\limits_{\varepsilon  \to \infty } &P \big\{ |\tilde r_c(t)-\tilde r_\iota | \le {C_0^{(\iota)}} {\rm e}^{-\gamma_0^{(\iota)}t} + \delta,\\
&\quad \forall t\in[0,T^{(\iota)}(\varepsilon)] \big\} =1 \text{\quad with} \mathop {\lim }_{\varepsilon  \to \infty }T^{(\iota)}(\varepsilon)=\infty,
\end{align*}
where ${C_0^{(\iota)}} = {c_0^{(\iota)}}| {{\bf\Xi} _\varepsilon ^{(\iota)}(0)}|.$ The proof is completed. \qed
\end{pf}
Theorem \ref{pro3} indicates the vehicle can locate a source in an elliptical signal map under small $V_c$. The vehicle finally points either directly towards or away from the source on the average. In fact the averaging heading of the vehicle is finally aligned with one of the coordinate axes. In other words, the vehicle converges to one point at the major or minor axis of the elliptical map. Note that the stability for the elliptical map under a constant forward velocity remains outstanding \cite{Cochran2009,Liu2010}, since in that case there is not a stable equilibrium to analyze (in the polar coordinates). For the case under large $V_c$, we cannot figure out an analytic solution due to the complexity of the average error system (\ref{eq:eplliticalEASys1})-(\ref{eq:eplliticalEASys4}), though simulation in Section \ref{sec:simulation} indicates the vehicle can also approach the source in this case. In Section \ref{subsec:circular} we shall study the result in a circular signal map, and compare it against the result in \cite{Liu2010}.

Next we give a brief discussion on the parameter selection for Theorem \ref{pro3}. Without loss of generality, we assume $q_p>0$. We also assume $a \in (0,3)$ and $g \in (0,3)$ to constrain the strength of the stochastic perturbation.

Using the fact that $2{\gamma _1}\left( {1 + {I_1}(2a,g)} \right) > {\gamma _2}{I_1}(a,g)$ and $q_r>2q_p$, one can easily derive that $V_c^l( \pm {q_p}) <  - \frac{{hR{\gamma _1}}}{{2I_1^2(a,g)}} < 0$. We also have $V_c^u( {\rm sgn}(-\gamma_3) {q_p}) > 0$ under the condition
\begin{equation}
b > \frac{{2cR}{\rm sgn}(\gamma_3){q_p}{\gamma _5}}{{({q_r} + 2{\mathop{\rm sgn}} (\gamma_3){q_p}){\gamma _4}}}. \label{eq:condition_v_c^u>0}
\end{equation}
Thus under condition (\ref{eq:condition_v_c^u>0}), we can always find an appropriate $V_c$ for Theorem \ref{pro3} by choosing $V_c$ small enough.

The sign of $V_c-V_{io}$ decides the average heading of the vehicle around the equilibria. The average heading would point inward when $V_c<V_{io}$ and outward when $V_c>V_{io}$. In addition, we have $V_c^l({q_p})<0<V_{io}$ when $\gamma_3<0$ and $V_c^l({-q_p})<V_{io}<0$ when $\gamma_3>0$. Observing that $\gamma_3>0$ if $a \in (0,1)$ and $\gamma_3<0$ if $a \in (1,3)$, we obtain the following corollary by summarizing the above analysis.
\begin{cor}
Consider the system in Theorem \ref{pro3} with $q_p>0$, $a \in (0,1) \cup (1,3)$ and $g \in (0,3)$, assume conditions in Theorem \ref{pro3} and (\ref{eq:condition_v_c^u>0}) are satisfied.
\par\noindent(i) When $a \in (0,1)$, the vehicle center converges to a point at the major axis of the elliptical level sets. Specially, the vehicle eventually points away from the source on the average under a small positive $V_c$.
\par\noindent(ii) When $a \in (1,3)$, the vehicle center converges to a point at the minor axis of the elliptical level sets. Specially, the vehicle eventually points towards the source on the average under a small negative $V_c$.
\end{cor}

\subsection{Stability for Circular Level Sets} \label{subsec:circular}
One can easily derive the stability for circular level sets under a small $V_c$ by setting $q_p=0$ in Theorem \ref{pro3}.
Due to the special structure of the circular level sets, we can also prove the stability under a large $V_c$. We write the signal distribution as
\begin{equation}
J = f(r_s) = {f^*} - {q_r}{\left| {{r_s} - {r^*}} \right|^2}, \label{eq:Jcircular}
\end{equation}
and rewrite the expression of $\xi$ as
\begin{equation*}
\xi =  - {q_r}\left( {{{\tilde r}_c}^2 - 2R{{\tilde r}_c}\cos (\hat \theta  - {\theta ^*} + a\chi(t/\varepsilon) )} \right) - {\tilde e_\xi }.
\end{equation*}
Observing the expression of $\xi$ and the shifted error system (\ref{eq:eplliticalSys1})-(\ref{eq:eplliticalSys5}), the system order can be reduced by defining $\tilde \theta  = \hat \theta  - {\theta ^*}$, which results in the following reduced shifted error system
\begin{align}
\frac{{\rm d}{{\tilde r}_c}}{{\rm d}t} &=  - ({V_c} + b\xi )\cos (\tilde \theta  + a\chi (\frac{t}{\varepsilon} )), \label{esys1}\\
\frac{{\rm d}\tilde \theta }{{\rm d}t} &= c\xi \sin (\chi (\frac{t}{\varepsilon} )) + \frac{{({V_c} + b\xi )}}{{{{\hat r}_c}}}\sin(\tilde \theta  + a\chi (\frac{t}{\varepsilon} )),\\
\frac{{\rm d}{{\tilde e}_\xi }}{{\rm d}t} &= h\xi,
\end{align}
\begin{align}
\xi  &=  - {q_r}{{\tilde r}_c}^2 - {{\tilde e}_\xi } + 2{q_r}R{{\tilde r}_c}\cos (\tilde \theta  + a\chi (\frac{t}{\varepsilon} )),\label{eq:xi}\\
{\rm d}\chi (t) &=  - \chi (t){\rm d}t + g{\rm d}B(t). \label{esys5}
\end{align}
The corresponding average error system is
\begin{align}
\frac{{{\rm d}\tilde r_c^{{\rm av}}}}{{{\rm d}t}} &= \left( {b{q_r}{{\tilde r}_c}^{{\rm av}2} + b\tilde e_\xi ^{{\rm av}} - {V_c}} \right)\cos ({{\tilde \theta }^{{\rm av}}}){I_1}(a,g) \notag\\
&\quad - b{q_r}R\tilde r_c^{{\rm av}}\cos (2{{\tilde \theta }^{{\rm av}}}){I_1}(2a,g) - b{q_r}R\tilde r_c^{{\rm av}}, \label{eq:av1}\\
\frac{{{\rm d}{{\tilde \theta }^{{\rm av}}}}}{{{\rm d}t}} &=  - 2c{q_r}R\tilde r_c^{{\rm av}}\sin({{\tilde \theta }^{{\rm av}}}){I_2}(a,g) \notag\\
&\quad + \frac{{{V_c} - b({q_r}{{\tilde r}_c}^{{\rm av}2} + \tilde e_\xi ^{{\rm av}})}}{{\tilde r_c^{{\rm av}}}}\sin({{\tilde \theta }^{{\rm av}}}){I_1}(a,g)\notag\\
&\quad + b{q_r}R\sin(2{{\tilde \theta }^{{\rm av}}}){I_1}(2a,g), \label{eq:av2}\\
\frac{{{\rm d}\tilde e_\xi ^{{\rm av}}}}{{{\rm d}t}} &=  - h{q_r}{{\tilde r}_c}^{{\rm av}2} - h\tilde e_\xi ^{{\rm av}} + 2h{q_r}R\tilde r_c^{{\rm av}}\cos ({{\tilde \theta }^{{\rm av}}}){I_1}(a,g), \label{eq:av3}
\end{align}
The average error system has four equilibria defined by\footnote{We have implicitly assumed ${{\tilde \theta }^{{\rm{ave}}}} \in (-\pi,\pi]$ to exclude repetitive equilibria.}
\begin{align}
\left[ {\tilde r_c^{{\rm{av^{eq1}}}},{{\tilde \theta }^{{\rm{av^{eq1}}}}},{{\tilde e_\xi}^{{\rm{av^{eq1}}}}}} \right] &= \left[ {\rho _1},\pi ,{e_1} \right], \label{eq:eq1}\\
\left[ {\tilde r_c^{{\rm{av^{eq2}}}},{{\tilde \theta }^{{\rm{av^{eq2}}}}},{{\tilde e_\xi}^{{\rm{av^{eq2}}}}}} \right] &= \left[ {-\rho _1},0,{e_1} \right], \label{eq:eq2}\\
\left[ {\tilde r_c^{{\rm{av^{eq3}}}},{{\tilde \theta }^{{\rm{av^{eq3}}}}},{{\tilde e_\xi}^{{\rm{av^{eq3}}}}}} \right] &= \left[ {{\rho _2},\alpha,{e_2}} \right], \label{eq:eq3}\\
\left[ {\tilde r_c^{{\rm{av^{eq4}}}},{{\tilde \theta }^{{\rm{av^{eq4}}}}},{{\tilde e_\xi}^{{\rm{av^{eq4}}}}}} \right] &= \left[ {{\rho_2},-\alpha,{e_2}} \right], \label{eq:eq4}
\end{align}
where
\begin{align*}
{\rho _1} &= \frac{{{V_c}{I_1}(a,g)}}{{b{q_r}R{\gamma _1}}},\\
{\rho _2} &= \sqrt {\frac{{c{V_c}{I_1}(a,g){I_2}(a,g) + {b^2}{q_r}R{\gamma _6}}}{{2{c^2}I_2^2(a,g){q_r}R}}},\\
{\alpha} &= \arccos \left( {\frac{{{\gamma _7}}}{{{\rho _2}}}} \right) \notag\\
&= \pi-\arctan\left( \frac{{\sqrt {{\gamma _8}} }}{{b\sqrt {{q_r}R} \left( {1 - {I_1}(2a,g)\;} \right)}}\right),\\
{e_1} &=  - \frac{{2{V_c}I_1^2(a,g)}}{{b{\gamma _1}}} - \frac{{V_c^2I_1^2(a,g)}}{{{q_r}{R^2}{b^2}\gamma _1^2}},\\
{e_2} &= 2{q_r}R{\gamma _7}{I_1}(a,g)- \frac{{c{V_c}{I_1}(a,g){I_2}(a,g) + {b^2}{q_r}R{\gamma _6}}}{{2{c^2}I_2^2(a,g)R}},\\
{\gamma _6} &= \left( {1 - {I_1}(2a,g)} \right)\left( {I_1^2(a,g) - {I_1}(2a,g)} \right),\\
{\gamma _7} &= \frac{{b\left( {{I_1}(2a,g) - 1} \right)}}{{2c{I_2}(a,g)}},\\
{\gamma _8} &= {2c{V_c}{I_1}(a,g){I_2}(a,g) - {b^2}{q_r}R{\gamma _4}}.\\ 
\end{align*}
Each of the equilibria (\ref{eq:eq1})-(\ref{eq:eq4}) represents an attractor around the source. Similar to the proof of Theorem \ref{pro3}, we can prove the local stability for the circular signal map.
\begin{thm} \label{pro2}
Consider the system (\ref{eq:sys1})-(\ref{eq:sys6}), (\ref{eq:Jcircular}) with positive parameters $a$, $g$, $b$, $c$, $h$, and $\varepsilon \in (0,\varepsilon_0)$.
\par\noindent (i) If the bias forward velocity $V_c$ is chosen such that
\begin{equation}
\text{either }  V_c \in (\bar V_c^l,0) \text{ or } V_c \in (0,\bar V_c^u), \label{cc:vc1}
\end{equation}
where
\begin{align*}
\bar V_c^l& \buildrel \Delta \over = V_c^l(0) = - \frac{{b{q_r}R{I_1}(2a,g) + b{q_r}R + h}}{{2I_1^2(a,g)}}R{\gamma _1}, \\
\bar V_c^u& \buildrel \Delta \over = V_c^u(0) =\frac{{{b^2}{q_r}R{\gamma _4}}}{{2c{I_1}(a,g){I_2}(a,g)}},
\end{align*}
then the system achieves the local convergence, both almost surely and in probability, to equilibrium (\ref{eq:eq1}) or equilibrium (\ref{eq:eq2}).
\par\noindent (ii) If the bias forward velocity $V_c$ is chosen such that
\begin{equation}
V_c > \bar V_c^u, \label{cc:vc3}
\end{equation}
then the system achieves the local convergence, both almost surely and in probability, to equilibrium (\ref{eq:eq3}) or equilibrium (\ref{eq:eq4}).
\end{thm}
\begin{pf}
The proof is given in Appendix \ref{appendix:proof}.
\end{pf}

Theorem \ref{pro2} describes the behavior of the vehicle in a circular signal map under different $V_c$.

When $V_c$ is small, the vehicle finally converges to an annular attractor near the source, whose radius can be small by choosing a small $V_c$ close to zero. In the meantime, the convergence rate would not decrease severely due to the tuning of the forward velocity. From (\ref{eq:eq1}) and (\ref{eq:eq2}) we infer that the limit of the vehicle's average heading is directly towards the source under small negative $V_c$ and away from the source under small positive $V_c$. The vehicle stays in the attractor as if it comes to a full stop. Note that in this case the vehicle may converge to any point at the annular attractor while in the elliptical case the attractor must be at the major or minor axis of the elliptical level sets. These phenomena are quite different from the result in \cite{Liu2010}, though both control scheme succeed in navigating the vehicle to the source with a partly random trajectory.

When $V_c$ is large, the vehicle converges to an annular attractor near the source while its average heading relative to the annulus is more outward than inward. It finally revolves around the source clockwise or counterclockwise on the average, depending on the initial conditions. This behavior is similar to the result in \cite{Liu2010}, but the vehicle revolves around the source with an non-tangential average heading due to the metabolic forward velocity.

\section{Simulation} \label{sec:simulation}
In this section we present simulation results to illustrate the behaviors of the vehicle in a signal map with circular or elliptical level sets. We also consider locating a source in a non-quadratic signal field. In all simulations we use band-limited white noise to approximate the white noise.

\subsection{Signal Maps with Circular Level Sets}
In this part, we examine the performance of the vehicle in a circular map. The map parameters are set as $f^* = 0$, $r^* = (0,0)$ and $q_r = 1.5$, and the initial conditions of the vehicle are set as $r_c(0) = (1,1)$, and $\theta(0) = -\pi/2$. The distance between the sensor and the vehicle center is set as $R = 0.1$. The controller parameters are chosen as $a = 2$, $g = 1$, $\varepsilon = 0.01$, $b = 2$, $c = 500$ and $h = 2$.

\begin{figure}[!t]
\centering
\subfigure[The trajectory of the vehicle center.]{\label{figEq1v}\includegraphics[width=0.8\hsize]{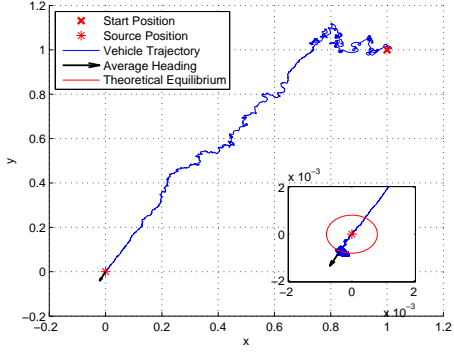}}
\subfigure[The sensor reading $J = f(r_s)$ and the forward velocity $v = V_c + b\xi$.]{\label{figEq1f}\includegraphics[width=0.8\hsize]{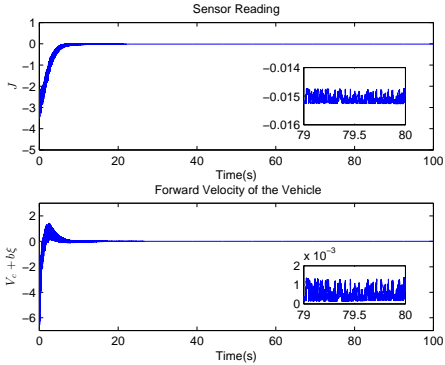}}
\caption{ Simulation results for circular level sets under small positive $V_c$.}
\label{figEq1}
\end{figure}

Fig.~\ref{figEq1} illustrates the behavior of the vehicle dictated by Theorem \ref{pro2} under small positive $V_c$. The bias forward velocity is chosen as $V_c = 0.0005$. As shown in Fig.~\ref{figEq1v}, the vehicle center converges to a small neighbourhood very close to the source with its heading points away from the source on the average. The trajectory of the vehicle center is partly random due to the using of the stochastic perturbation. Fig.~\ref{figEq1f} shows the sensor reading and the forward velocity of the vehicle.

%

In Fig.~\ref{figEq1} the vehicle finally moves in a small area near the source as if it comes to a full stop. This is quite different from the result of \cite{Liu2010}, where the vehicle finally drifts in an annulus around the source. In addition, the attractor in Fig.~\ref{figEq1} is very close to the source under a small $V_c$ since ${{\tilde r_c }^{{\rm{av^{eq1}}}}}$ is positively correlated to $|V_c|$. Note that the vehicle does not strictly stop at the source, though its stop seems evident from Fig.~\ref{figEq1v}. We can see it keeps moving in the attractor with a small forward velocity from Fig.~\ref{figEq1f}.

\begin{figure}[!t]
  \centering
  \includegraphics[width=0.8\hsize]{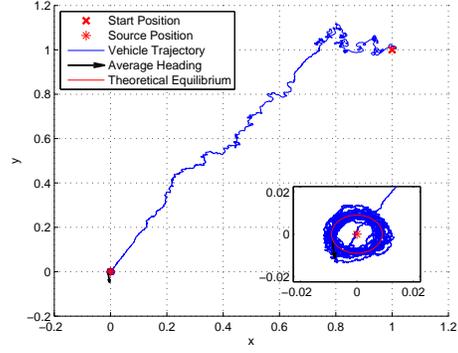}
  \caption{Vehicle trajectory for circular level sets under large $V_c$.}
  \label{figEq3}
\end{figure}

Fig.~\ref{figEq3} illustrates the behavior of the vehicle dictated by Theorem \ref{pro2} under large $V_c$. The bias forward velocity is chosen as $V_c = 0.01$. The vehicle converges to an annular attractor and revolves around the source, which is similar to the result in \cite{Liu2010}. The average heading is more outward than inward, which coincides with the theoretical result.

\subsection{Signal Maps with Elliptical Level Sets}
In this part, we examine the performance of the vehicle in an elliptical map. The map parameters are set as $f^* = 0$, $r^* = (0,0)$, $q_r = 2$ and $q_p = 0.5$ and the initial conditions of the vehicle are set as $r_c(0) = (1,1)$, and $\theta(0) = -\pi/2$. The distance between the sensor and the vehicle center is set as $R = 0.1$. The controller parameters are chosen as $\varepsilon = 0.01$, $b = 2$, $c = 500$ and $h = 2$.

\begin{figure}[!t]
  \centering
  \includegraphics[width=0.8\hsize]{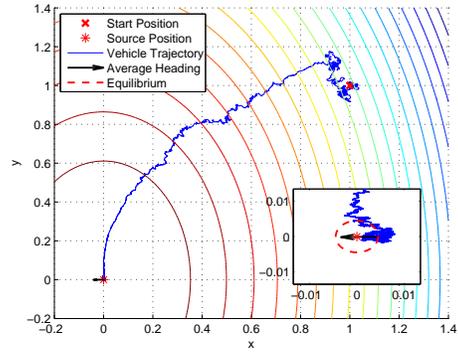}
  \caption{Vehicle trajectory for elliptical level sets under small $V_c$.}
  \label{figElliptical}
\end{figure}

\begin{figure}[!t]
  \centering
  \includegraphics[width=8.5cm]{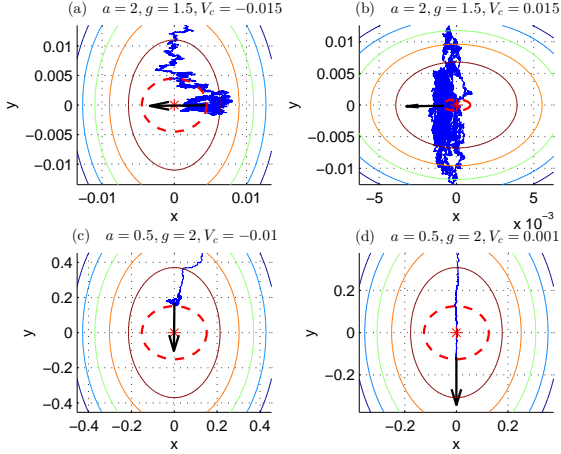}
  \caption{Vehicle trajectories for elliptical level sets under different parameters.}
  \label{figEllipticalAll}
\end{figure}

Fig.~\ref{figElliptical} and Fig.~\ref{figEllipticalAll} illustrate the behavior of the vehicle dictated by Theorem \ref{pro3}. In Fig.~\ref{figElliptical} and Fig.~\ref{figEllipticalAll}(a), we chose $a = 2$, $g = 1.5$ and $V_c = -0.015$. In Fig.~\ref{figEllipticalAll}(b), we chose $a = 2$, $g = 1.5$ and $V_c = 0.015$. In Fig.~\ref{figEllipticalAll}(c), we chose $a = 0.5$, $g = 2$ and $V_c = -0.01$. In Fig.~\ref{figEllipticalAll}(d), we chose $a = 0.5$, $g = 2$ and $V_c = 0.001$. As depicted in Fig.~\ref{figElliptical} and Fig.~\ref{figEllipticalAll}, the vehicle converges to a small area near the source under small $V_c$. Fig.~\ref{figEllipticalAll} illustrates the convergence to different equilibria under different parameters in the same signal map.

\begin{figure}[!t]
  \centering
  \includegraphics[width=0.8\hsize]{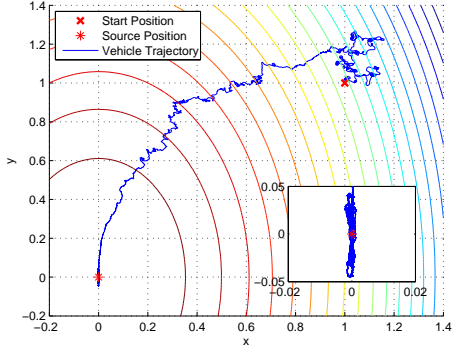}
  \caption{Vehicle trajectory for elliptical level sets under large $V_c$.}
  \label{figEllipticalLarge}
\end{figure}

Fig.~\ref{figEllipticalLarge} illustrates the behavior of the vehicle under large $V_c$. The parameters are chosen as $a = 2$, $g = 1$, $q_p=0.5$ and $V_c=0.01$. In this case, the vehicle can also approaches the source. As shown in Fig.~\ref{figEllipticalLarge}, it overshoots the source, turns back and overshoots the source again, and so on.

\subsection{Non-Quadratic Signal Maps}
\begin{figure}[!t]
  \centering
  \includegraphics[width=0.8\hsize]{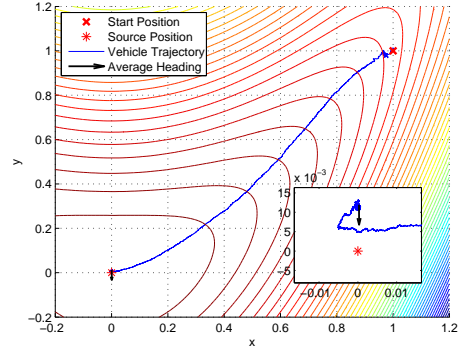}
  \caption{Vehicle trajectory for a Rosenbrock function signal map.}
  \label{figRosenbrock}
\end{figure}

Our control scheme also exhibits abilities to seek the sources of signal fields with non-quadratic maps. In Fig.~\ref{figRosenbrock} we assume the signal distribution is a Rosenbrock function, which takes the form $J = -x_s^2-(y_s-x_s^2)^2$. The Rosenbrock function has an isolated maximum at $(0,0)$ and its Hessian at $(0,0)$ is negative definite. The initial conditions of the vehicle and controller parameters are chosen to be the same as those in Fig.~\ref{figEq1} except $V_c = -0.0005$. As depicted in Fig.~\ref{figRosenbrock}, the vehicle can also well approach the source.

\section{Conclusion}
We have studied the nonholonomic source seeking problem in a plane. In our control scheme, both forward and angular velocities are tuned according to the stochastic extremum method\cite{liu2010stochastic}. As a result, the vehicle well approaches the source with a partly random trajectory. We adopted the stochastic averaging theory for nonlinear continuous-time systems to prove the local stability for static signal fields with elliptical level sets. We have established the local exponential convergence, both almost surely and in probability, to attractors in an annulus around the source. Under a small bias forward velocity the vehicle may virtually ``stop" at the source without sacrificing the convergence rate.


\appendix
\section{Proof of Theorem \ref{pro2}} \label{appendix:proof}
We complete the proof of the first part of Theorem \ref{pro2} in Section \ref{appendix:proof1}, which establishes the convergence to equilibria (\ref{eq:eq1}) and (\ref{eq:eq2}). Then we prove the second part in Section \ref{appendix:proof2}, which corresponds to equilibria (\ref{eq:eq3}) and (\ref{eq:eq4}).
\subsection{Under Small $V_c$}\label{appendix:proof1}
The Jacobians of equilibria (\ref{eq:eq1}) and (\ref{eq:eq2}) are as follows, respectively,
\begin{align*}
{A}^{\rm eq1} = \begin{bmatrix}
{m_{11}}&0&{ -m_{13}}\\
0&{m_{22}}&0\\
{-m_{31}}&0&{-h}
\end{bmatrix},
{A}^{\rm eq2} = \begin{bmatrix}
{{m_{11}}}&0&{{m_{13}}}\\
0&{{m_{22}}}&0\\
{{m_{31}}}&0&{-h}
\end{bmatrix},
\end{align*}
where
\begin{align*}
{m_{11}} &=  - 2b{q_r}{I_1}(a,g){\rho _1}\; - b{q_r}R{I_1}(2a,g) - b{q_r}R,\\
{m_{13}} &= b{I_1}(a,g),\\
{m_{22}} &= \frac{{2c{V_c}{I_1}(a,g){I_2}(a,g)}}{{b{\gamma _1}}} - b{q_r}R\left( {1 - {I_1}(2a,g)} \right),\\
{m_{31}} &= 2h{q_r}{\rho _1} + 2h{q_r}R{I_1}(a,g).
\end{align*}

The two Jacobians have the same characteristic equation, which is given by
\begin{equation}
\left( {\lambda  - {m_{22}}} \right)\left( {{\lambda ^2} + (h-{m_{11}})\lambda  - {m_{11}}h - {m_{13}}{m_{31}}} \right) = 0. \label{eq:characteristic}
\end{equation}

To guarantee that all roots of characteristic equation (\ref{eq:characteristic}) have negative real parts, we need
\begin{align}
{m_{22}} < 0, \label{eq:cc1}\\
{m_{11}} - h < 0, \label{eq:cc2}\\
-{m_{11}}h - {m_{13}}{m_{31}} > 0. \label{eq:cc3}
\end{align}
The inequality (\ref{eq:cc3}) holds for all $\gamma_1>0$ and the inequalities (\ref{eq:cc1}) and (\ref{eq:cc2}) are satisfied under condition (\ref{cc:vc1}). Hence, the Jacobians ${A}^{\rm eq1}$ and ${A}^{\rm eq2}$ are Hurwitz under condition (\ref{cc:vc1}), which implies
that both equilibria (\ref{eq:eq1}) and (\ref{eq:eq2}) are exponentially stable. Applying Theorem 2 in \cite{liu2010stochastic}, we finish the proof through an inference similar to the latter part of the proof of Theorem \ref{pro3}.
\subsection{Under Large $V_c$}  \label{appendix:proof2}
We first prove that under condition (\ref{cc:vc3}) equilibria (\ref{eq:eq3}) and (\ref{eq:eq4}) are valid, i.e., $\gamma_{\rho_2} \buildrel \Delta \over = c{V_c}{I_1}(a,g){I_2}(a,g) + {b^2}{q_r}R{\gamma_6} > 0$ and $\gamma_8>0$. When $V_c > \bar V_c^u$, we have
\begin{align*}
\gamma_{\rho_2}& > c{\bar V_c^u}{I_1}(a,g){I_2}(a,g) + {b^2}{q_r}R{\gamma _6}\\
& = \frac{1}{2}{b^2}{q_r}R{\left( {1 - {I_1}(2a,g)} \right)^2}\quad >0,\\
{\gamma_8} &> 2c{\bar V_c^u}{I_1}(a,g){I_2}(a,g) - {b^2}{q_r}R{\gamma _4} = 0.
\end{align*}
Thus the two equilibria are well defined.

The Jacobians of equilibria (\ref{eq:eq3}) and (\ref{eq:eq4}) are given by
\begin{align*}
{A}^{\rm eq3} &= \begin{bmatrix}
{{k_{11}}}&{{k_{12}}}&{{k_{13}}}\\
{{k_{21}}}&{{k_{22}}}&{{k_{23}}}\\
{{k_{31}}}&{{k_{32}}}&{ - h}
\end{bmatrix},
{A}^{\rm eq4} &= \begin{bmatrix}
{{k_{11}}}&{ - {k_{12}}}&{{k_{13}}}\\
{ - {k_{21}}}&{{k_{22}}}&{ - {k_{23}}}\\
{{k_{31}}}&{ - {k_{32}}}&{ - h}
\end{bmatrix},
\end{align*}
where the coefficients of the Jacobians are provided in Appendix \ref{appendix:coefficient}. The two Jacobians have the same characteristic equation, which is given by
\begin{equation}
{\lambda ^3} + {l_2}{\lambda ^2} + {l_1}\lambda  + l_0 = 0, \label{eq:characteristic2}
\end{equation}
where
\begin{align*}
{l_0} &= 2h{q_r}R{\gamma _8},\\
{l_1} &= 2{q_r}R{\gamma _8} + b{q_r}\frac{{{I_1}(a,g)}}{{c{I_2}(a,g)}}{\gamma _8} + hb{q_r}R{\gamma _1},\\
{l_2} &= h + b{q_r}R\left( {{I_1}(2a,g) + 1} \right) - 2b{q_r}{\gamma _7}{I_1}(a,g).
\end{align*}
Invoking the Routh-Hurwitz test, it's easily to verify the Jacobians ${A}^{\rm eq2}$ and ${A}^{\rm eq3}$ are Hurwitz under condition (\ref{cc:vc3}). Hence, equilibria (\ref{eq:eq1}) and (\ref{eq:eq2}) are exponentially stable. Applying Theorem 2 in \cite{liu2010stochastic}, we finish the proof through an inference similar to the latter part of the proof of Theorem \ref{pro3}.
\section{Coefficients of the Jacobians} \label{appendix:coefficient}
The coefficients of $J_1(q_p)$ and $J_2(q_p)$ are as follows
\begin{align*}
{a_{11}}({q_p}) &= 2b\left( {{q_r} + 2{q_p}} \right)\rho ({q_p}){I_1}(a,g) \\
&\quad - bR\left( {{q_r} + 2{q_p}} \right)\left( {1 + {I_1}(2a,g)} \right),\\
{a_{14}} &= b{I_1}(a,g),\\
{a_{22}}({q_p}) &=  - bR({q_r} - 2{q_p})(1 - {I_1}(2a,g)),\\
{a_{23}}({q_p}) &= 2bR({q_r} + 2{q_p})\left( {I_1^2(a,g) - {I_1}(2a,g)} \right) + \frac{1}{{\rho ({q_p})}}\\
&\quad \times \Big( ({3{I_1}(3a,g) - {I_1}(a,g) - 2{I_1}(2a,g){I_1}(a,g)} ) \\
&\quad\quad\quad \times b{q_p}{R^2} - {V_c}{I_1}(a,g) \Big),\\
{a_{32}}({q_p}) &= 2cR\left( {{q_r} - 2{q_p}} \right)\rho ({q_p}){I_2}(a,g),\\
{a_{33}}({q_p}) &= 4c{q_p}{R^2}{I_2}(2a,g) - 2cR\left( {{q_r} + 2{q_p}} \right)\rho ({q_p}){I_2}(a,g),\\
{a_{41}}({q_p}) &= 2hR\left( {{q_r} + 2{q_p}} \right){I_1}(a,g) - 2h\left( {{q_r} + 2{q_p}} \right)\rho ({q_p}).
\end{align*}
The coefficients of ${A}^{\rm eq3}$ and ${A}^{\rm eq4}$ are as follows
\begin{align*}
{k_{11}} &= 2b{q_r}{\gamma _7}{I_1}(a,g) - 2b{q_r}R{\cos ^2}({\alpha}){I_1}(2a,g) \\
&\quad + b{q_r}R\left( {{I_1}(2a,g) - 1} \right),\\
{k_{12}} &= 2c{q_r}R\rho _2^2\sin({\alpha}){I_2}(a,g) + b{q_r}R{\rho _2}\sin(2{\alpha}){I_1}(2a,g)\\
{k_{13}} &= b\cos ({\alpha}){I_1}(a,g),\\
{k_{21}} &=  - 4c{q_r}R\sin({\alpha}){I_2}(a,g) - 2b{q_r}\sin({\alpha}){I_1}(a,g)\\
&\quad + b{q_r}R\frac{{\sin(2{\alpha})}}{{{\rho _2}}}{I_1}(2a,g),\\
{k_{22}} &=  - 2b{q_r}R{\sin ^2}({\alpha}){I_1}(2a,g),\\
{k_{23}} &= \frac{{ - b}}{{{\rho _2}}}\sin({\alpha}){I_1}(a,g),\\
{k_{31}} &=  - 2h{q_r}{\rho _2} + 2h{q_r}R\cos ({\alpha}){I_1}(a,g),\\
{k_{32}} &=  - 2h{q_r}R{\rho _2}\sin({\alpha}){I_1}(a,g).
\end{align*}

\bibliographystyle{plain}
\bibliography{myAutoRef}

\end{document}